\documentclass{elsart3-1}

\usepackage{amssymb}

\usepackage[english,francais]{babel}

\newtheorem{theorem}{Theorem}[section]
\newtheorem{lemma}[theorem]{Lemma}
\newtheorem{e-proposition}[theorem]{Proposition}

\newtheorem{e-definition}[theorem]{Definition\rm}

\renewcommand{\theequation}{\arabic{equation}}
\def\og{\leavevmode\raise.3ex\hbox{$\scriptscriptstyle\langle\!\langle$~}}
\def\fg{\leavevmode\raise.3ex\hbox{~$\!\scriptscriptstyle\,\rangle\!\rangle$}}

\journal{the Acad\'emie des sciences}
\begin{document}
\centerline{}
\begin{frontmatter}


\selectlanguage{english}
\title{A Berry-Esseen bound of order $ 1/\sqrt{n} $ for martingales}


\selectlanguage{english}

\author[label1]{Songqi Wu}
\author[label1]{Xiaohui Ma}
\author[label2]{Hailin Sang}
\author[label1]{Xiequan Fan}
\ead{fanxiequan@hotmail.com}
\address[label1]{Center for Applied Mathematics,
Tianjin University, Tianjin, China}
\address[label2]{Department of Mathematics, The University of Mississippi, University, MS 38677, USA}


\medskip
\begin{center}
{\small Received *****; accepted after revision +++++\\
Presented by £££££}
\end{center}

\begin{abstract}
\selectlanguage{english}
Renz \cite{R96} has established a rate of convergence $1/\sqrt{n}$ in the central limit theorem for martingales with some restrictive conditions. In the present paper a modification of the methods, developed by Bolthausen \cite{B82} and Grama and Haeusler \cite{GH}, is applied for obtaining the same convergence rate for a class of more general martingales. An application to linear processes is discussed.
{\it To cite this article: A.
Name1, A. Name2, C. R. Acad. Sci. Paris, Ser. I 340 (2005).}

\vskip 0.5\baselineskip

\selectlanguage{francais}
\noindent{\bf R\'esum\'e} \vskip 0.5\baselineskip \noindent
{\bf Une limite de Berry-Esseen de l$'$ordre $ 1 / \sqrt{n} $ pour les martingales.}
Renz  \cite{R96} a \'{e}tabli un taux de convergence $1/\sqrt{n}$ dans le th\'{e}or\`{e}me de la limite centrale pour les martingales avec certaines conditions restrictives. Dans le pr\'{e}sent article, une modification des m\'{e}thodes, d\'{e}velopp\'{e}es par Bolthausen  \cite{B82}  et Grama et Haeusler \cite{GH}, a appliqu\'{e} pour obtenir le m\^{e}me taux de convergence pour une classe de martingales plus g\'{e}n\'{e}rales.
Une application aux processus lin\'{e}aires est discut\'{e}e.
{\it Pour citer cet article~: A. Name1, A. Name2, C. R. Acad. Sci.
Paris, Ser. I 340 (2005).}
\end{abstract}
\end{frontmatter}

\section{Introduction and main result}
For $ n \in \mathbf{N} $, let $(\xi_i,\mathcal{F}_i)_{i=0,...,n}$ be a finite sequence of martingale differences defined on some probability space $(\Omega,\mathcal{F},\mathbf{P})$, where $\xi_0 = 0$ and $\{\emptyset,\Omega\} = \mathcal{F}_0 \subseteq...\subseteq \mathcal{F}_n \subseteq \mathcal{F}$ are increasing $\sigma$-fields. Denote \[X_0 = 0,\quad X_k = \sum_{i=1}^{k}\xi_i,\ k = 1,...,n.\]
Then $X = (X_k,\mathcal{F}_k)_{k=0,...,n}$ is a martingale. Denote by $\langle X \rangle$   the conditional variance of $X$:\[\langle X \rangle_0 = 0, \quad \langle X \rangle_k = \sum_{i=1}^{k}\mathbf{E}\big[\xi_i^2\big|\mathcal{F}_{i-1}\big],\ k = 1,...,n.\]
Define \[D(X_n) = \sup_{x \in \mathbf{R}}\Big|\mathbf{P}(X_n \le x) - \Phi(x)\Big|,\]
where $\Phi(x)$ is the distribution function of the standard normal random variable. Denote by $\stackrel{\mathbf{P}}{\rightarrow}$ the convergence in probability as $n \rightarrow \infty$. According to the martingale central limit theorem, the ``conditional Lindeberg condition" \[
\sum_{i=1}^{n}\mathbf{E}\big[\xi_i^2\mathbf{1}_{\{|\xi_i| \ge \epsilon\}}\big|\mathcal{F}_{i-1}\big] \stackrel{\mathbf{P}}{\rightarrow} 0,\quad \textrm{for each} \ \epsilon > 0,\]
and the ``conditional normalizing condition"
$\langle X \rangle_n \stackrel{\mathbf{P}}{\rightarrow}1$
together implies asymptotic normality of $ X_n, $ that is,
$D(X_n) \rightarrow0$ as $n \rightarrow \infty.$


The convergence rate of $ D(X_n) $ has attracted a lot of attentions. For instance, Bolthausen \cite{B82} proved that if $ |\xi_i| \le \epsilon_n $ for a number $ \epsilon_n $ and $ \langle X \rangle_n = 1 $ a.s., then $ D(X_n) \le c\epsilon_n^3n\log n, $ where, here and after, $ c $ is an absolute constant not depending on $ \epsilon_n $ and $ n $. El Machkouri and Ouchti \cite{M07} improved the factor $ \epsilon_n^3n\log n $ in Bolthausen's bound to $ \epsilon_n\log n $ under the following more general condition \[ \mathbf{E}\big[|\xi_i|^3\big|\mathcal{F}_{i-1}\big] \le \epsilon_n\mathbf{E}\big[\xi_i^2\big|\mathcal{F}_{i-1}\big] \quad  a.s.\ \textrm{for all} \ i = 1,2,...,n.\]For more related results, we refer to Ouchti \cite{O05} and Mourrat \cite{M13}. Recently, Fan \cite{F19} proved that if there exist a positive constant $\rho$ and a number $\epsilon_n,$ such that \[ \mathbf{E}\big[|\xi_i|^{2+\rho}\big|\mathcal{F}_{i-1}\big] \le \epsilon_n^\rho\mathbf{E}\big[\xi_i^2\big|\mathcal{F}_{i-1}\big] \quad a.s.\ \textrm{for all}\ i = 1,2,...,n, \]  and $\langle X \rangle_n = 1 $ a.s., then $ D(X_n) \le c_\rho\hat{\epsilon}_n, $ where
\begin{displaymath}
 \hat{\epsilon}_n =
 \left\{ \begin{array}{ll}
 \epsilon_n^\rho, & \textrm{if $\rho \in (0,1)$}\\
 \epsilon_n|\log \epsilon_n|, & \textrm{if $\rho \geq 1$}
\end{array} \right.
\end{displaymath}
and $c_\rho$ is a constant depending only on $\rho$. Fan also showed that this Berry-Esseen bound is optimal. In particular, if $ \epsilon_n \asymp 1/\sqrt{n}, $ then we have $ \epsilon_n|\log \epsilon_n| \asymp (\log n)/\sqrt{n}. $ Thus, we cannot obtain the classical convergence rate $ 1/\sqrt{n} $ for general martingales.

However, the convergence rate $ 1/\sqrt{n} $ for martingales is possible to be attained with some additional restrictive conditions.  For instance,  Renz \cite{R96} proved that if there exists a constant $ \rho > 0 $ such that
\begin{equation}\label{Renz}
	\mathbf{E}[\xi_i^2|\mathcal{F}_{i-1}] = 1/n ,\ \ \ \ \mathbf{E}[\xi_i^3|\mathcal{F}_{i-1}] = 0 \ \ \ \  \textrm{and} \ \ \ \ \mathbf{E}\big[|\xi_i|^{3+\rho}\big|\mathcal{F}_{i-1}\big] \le c{n^{-(3+\rho)/2}}, \quad \textrm{a.s.},
\end{equation}
then it holds
\begin{equation}\label{Renzbound}
	D(X_n) = O\bigg(\frac{1}{\sqrt{n}}\bigg).
\end{equation}
He also showed that this result is not true for $ \rho = 0 .$
More martingale Berry-Esseen bounds   of convergence rate $ 1/\sqrt{n}$ can  be found in Bolthausen \cite{B82} and Kir'yanova and Rotar \cite{KR91}.

In this paper we are interested in extending (\ref{Renzbound}) to a class of more general martingales.
The following theorem is our main result.

\begin{theorem}\label{th1}
Assume that there exist some numbers $\rho\in(0,+\infty), \epsilon_n\in (0,\frac{1}{2}]$ and $\delta_n \in [0,\frac{1}{2}]$ such that for all $1 \le i \le n,$
\begin{eqnarray}\label{condi1}
\big|\langle X \rangle_n - 1\big|\le\delta_n^2,  	
\end{eqnarray}
\begin{eqnarray}\label{condi2}
\mathbf{E}\big[\xi_i^3\big|\mathcal{F}_{i-1}\big] = 0
\end{eqnarray}
and
\begin{eqnarray}\label{condi3}
\mathbf{E}\big[|\xi_i|^{3+\rho}\big|\mathcal{F}_{i-1}\big] \le \epsilon_n^{1+\rho}\mathbf{E}\big[\xi_i^2\big|\mathcal{F}_{i-1}\big] \quad \textrm{a.s. }
\end{eqnarray}
Then \[D(X_n)\le c_{\rho}(\epsilon_n + \delta_n),\] where $ c_\rho $ depends only on $ \rho. $ In addition, it holds $ c_{\rho} = O(\rho^{-1}), \rho \rightarrow 0. $
\end{theorem}

Notice that under the conditions of Renz \cite{R96}, the conditions of Theorem \ref{th1} are satisfied with $ \delta_n = 0 $ and $ \epsilon_n \asymp 1/\sqrt{n} .$ Thus Theorem \ref{th1} extends Renz's result to a class of more general martingales.

Thanks to the additional condition (\ref{condi2}), the Berry-Esseen bound (\ref{5}) improves the bound of Fan \cite{F19} by replacing $ \epsilon_n|\log \epsilon_n| $ with $ \epsilon_n $.

Relaxing the condition  (\ref{condi1}), we have the following analogue estimation of Fan (cf.\ (26) of \cite{F19}).
\begin{theorem}\label{th2}
 Assume that there exist some numbers $\rho\in(0,+\infty)$ and $\epsilon_n\in(0,\frac{1}{2}]$ such that for all $1\leq i\leq n,$
\[\mathbf{E}\big[\xi_i^3\big|\mathcal{F}_{i-1}\big]=0\] and
\[\mathbf{E}\big[|\xi_i|^{3+\rho}\big|\mathcal{F}_{i-1}\big]\leq \epsilon_n^{1+\rho}\mathbf{E}\big[\xi_i^2\big|\mathcal{F}_{i-1}\big] \quad a.s.\]
Then, for all $p\geq 1$,
\begin{equation}\label{5}
D(X_n)\leq c_{\rho} {\epsilon_n}+c_{p}\bigg(\mathbf{E}\Big[\big|\langle X\rangle_n - 1\big|^p\Big]+\mathbf{E}\Big[\max_{1\leq i\leq n}|\xi_i|^{2p}\Big]\bigg)^{1/{(2p+1)}} ,
\end{equation}
where $c_{\rho}$ and $c_{p}$ depend only on $\rho$ and $p$, respectively.
\end{theorem}

It is easy to see that when $p \rightarrow \infty,$
\[ \bigg(\mathbf{E}\Big[\big|\langle X\rangle_n - 1\big|^p\Big]\bigg)^{1/{(2p+1)}} \rightarrow ||\langle X \rangle_n - 1 ||_{\infty}^{1/2}, \]
which coincides with $ \delta_n $ of Theorem \ref{th1}.


\section{Application}
We first extend Theorem 1.1 to triangular arrays with infinity many terms in each line. For $n\in \mathbf{N}$, let $(\xi_{n,i}, \mathcal{F}_{n,i})_{i=-\infty}^n$ be a sequence of martingale differences defined on some probability space $(\Omega, \mathcal{F}, \mathbf{P})$, where the adapted filtration is $\{\varnothing, \Omega\}=\mathcal{F}_{-\infty}\subset ...\subset \mathcal{F}_{n, n-1}\subset \mathcal{F}_{n, n}\subset \mathcal{F}$. Denote $X_{n,k}=\sum_{i=-\infty}^k\xi_{n,i}, k\le n$. Then $(X_{n,k}, \mathcal{F}_{n,k})_{k=-\infty}^n$ is a martingale. Let $\langle X\rangle_{n,k}=\sum_{i=-\infty}^k \mathbf{E}[\xi_{n,i}^2|\mathcal{F}_{n, i-1}], k\le n$.  In particular, denote  $X_n:=X_{n,n}$ and $\langle X\rangle_{n}:=\langle X\rangle_{n,n}$.

With some slight modification on the proof, Theorem 1.1 still holds in this new setting. Now we apply Theorem 1.1 with this new setting to the partial sum of linear processes. Let $(\varepsilon_i)_{i\in\mathbf{Z}}$ be a sequence of identically distributed martingale differences adapted to the filtration $(\mathcal{F}_i)_{i\in\mathbf{Z}}$. We consider the causal linear process in the form
\begin{equation}
Y_{k}=\sum_{j=-\infty}^{k}a_{k-j}\varepsilon_{j}, \label{defx}%
\end{equation}
where the martingale differences have finite variance and the sequence of real coefficients satisfies $\sum_{i=0}^\infty a_i^2<\infty$. Without loss of generality, let the variance of the martingale difference to be $1$. We say the linear process has long memory if $\sum_{i=0}^\infty|a_i|=\infty$. In this case, we assume that $a_0=1$ and
\begin{equation}
a_{i}=l(i)i^{-\alpha},\textrm{ }i>0,\textrm{ with }1/2<\alpha<1. \label{rv}%
\end{equation}
Here $l(\cdot)$ is a slowly varying function. On the other hand, we say the linear process has short memory if $\sum_{i=0}^\infty |a_i|<\infty$ and $\sum_{i=0}^\infty a_i\ne 0$. The third case is $\sum_{i=0}^\infty |a_i|<\infty$ and $\sum_{i=0}^\infty a_i= 0$.

The long memory linear processes covers the well-known
fractional ARIMA processes (cf. Granger and Joyeux \cite{GrangerJoyeux}; Hosking \cite{Hosking}),
which play an important role in financial time series modeling and
application. As a special case, let $0<d<1/2$ and $B$ be the backward shift
operator with $B\varepsilon_{k}=\varepsilon_{k-1}$ and consider
\[
Y_{k}=(1-B)^{-d}\varepsilon_{k}=\sum_{i=0}^\infty a_{i}\varepsilon_{k-i},\textrm{ where }a_{i}%
=\frac{\Gamma(i+d)}{\Gamma(d)\Gamma(i+1)}.
\]
For this example we have $\lim_{n\rightarrow\infty}a_{n}/n^{d-1}=1/\Gamma(d)$.
Note that these processes have long memory because $\sum_{j=0}^\infty
|a_{j}|=\infty$.

The partial sum $S_n=\sum_{k=1}^n Y_k$ of causal linear process (\ref{defx})  can be written as
$S_n=\sum_{-\infty}^{n}b_{n,i}\varepsilon_i$, where $b_{n,i}=\sum_{j=0}^{n-i}a_j$ for $0<i\le n$, and  $b_{n,i}=\sum_{j=1-i}^{n-i}a_j$ for $i\le 0$. The variance of $S_n$ is $B_n^2=var(S_n)=\sum_{-\infty}^{n} b_{n,i}^2$.  Now let $X_{n,k}=\sum_{-\infty}^{k} b_{n,i}\varepsilon_i/B_n$. Then $X_n=X_{n,n}=S_n/B_n$ and $\langle X\rangle_n=\sum_{i=-\infty}^n b_{n,i}^2\mathbf{E}[\varepsilon_{i}^2|\mathcal{F}_{i-1}]/B_n^2$. If we assume $|\langle X\rangle_n-1|\le \delta_n^2$ for some $\delta_n\in [0,\frac{1}{2}]$, $\mathbf{E}[\varepsilon_i^3|\mathcal{F}_{i-1}]=0$ and $\mathbf{E}[|\varepsilon_i|^{3+\rho}|\mathcal{F}_{i-1}]\le d_\rho^{1+\rho}\mathbf{E}[\varepsilon_i^2|\mathcal{F}_{i-1}]$ a.s.\ for all $i\in \mathbf{Z}$ and some constant $d_\rho$, then, by Theorem 1.1, $$\sup_{x\in \mathbf{R}}| \mathbf{P}(S_n/B_n\le x)-\Phi(x)|\le c_\rho(\epsilon_n+\delta_n),$$ where $\epsilon_n=d_\rho\sup_{i\le n}|b_{n,i}|/B_n$.

In the case that $\sum_{i=0}^\infty |a_i|<\infty$, $\sup_{i\le n}|b_{n,i}|\le \sum_{i=0}^\infty |a_i|<\infty$ and it is well known that $B_n^2$ has order $n$. Hence $\epsilon_n$ has order $1/\sqrt{n}$ in this case. In the long memory case $\sum_{i=0}^\infty |a_i|=\infty$, if we assume (\ref{rv}), $B_n^2$ has order $n^{3-2\alpha}l^2(n)$ (e.g., Wu and Min \cite{WM}) and $\sup_{i\le n}|b_{n,i}|$ has order $n^{1-\alpha}l(n)$ (see Beknazaryan et al. \cite{BSX} for upper bound and Fortune et al. \cite{FPS} for lower bound in the case $d=1$). Hence in this case $\epsilon_n$ also has order $1/\sqrt{n}$.  In either case the Berry-Esseen bound has order $1/\sqrt{n}$ if $\delta_n=O(n^{-1/2})$. In particular if we in addition assume that the innovations $(\varepsilon_i)_{i\in\mathbf{Z}}$ are independent, then $\delta_n=0$ and the Berry-Esseen bound $\sup_{x\in \mathbf{R}}| \mathbf{P}(S_n/B_n\le x)-\Phi(x)|$ has order $1/\sqrt{n}$. Here the condition $\mathbf{E}[\varepsilon_i^3|\mathcal{F}_{i-1}]=0$ is needed to have  the Berry-Esseen bound of order $1/\sqrt{n}$. 
We cannot have this order from the result  of   Fan \cite{F19}.

\section{Proofs of theorems}
\subsection{Preliminary lemmas}\label{sec2}

In the proofs of theorems, we need the following technical lemmas. The first two lemmas can be found in Fan \cite{F19} (cf. Lemmas 3.1 and 3.2 therein).   
\begin{lemma} \label{lemma1}
	If there exists an $s > 3$  such that \[\mathbf{E}[|\xi_i|^s|\mathcal{F}_{i-1}] \le \epsilon_n^{s-2}\mathbf{E}[\xi_i^2|\mathcal{F}_{i-1}],\]
	then, for any $t \in [3,s)$,
	\[\mathbf{E}[|\xi_i|^t|\mathcal{F}_{i-1}] \le \epsilon_n^{t-2}\mathbf{E}[\xi_i^2|\mathcal{F}_{i-1}].\]
\end{lemma}
\begin{lemma}\label{lemma2}
	If there exists an $s>3$  such that \[\mathbf{E}[|\xi_i|^s|\mathcal{F}_{i-1}] \le \epsilon_n^{s-2}\mathbf{E}[\xi_i^2|\mathcal{F}_{i-1}],\]
	then \[\mathbf{E}[\xi_i^2|\mathcal{F}_{i-1}] \le \epsilon_n^2.\]
\end{lemma}
The next two technical lemmas are due to Bolthausen (cf. Lemmas 1 and 2 of \cite{B82}).
\begin{lemma}\label{lemma3}
	Let $X$ and $Y$ be random variables. Then \[\sup_u\Big|\mathbf{P}\big(X \le u\big) - \Phi(u)\Big| \le c_1\sup_u\Big|\mathbf{P}\big(X + Y \le u\big) - \Phi(u)\Big| + c_2\Big|\Big|\mathbf{E}[Y^2|X]\Big|\Big|_{\infty}^{1/2},\]
	where $ c_1$ and $ c_2 $ are two positive constants.	
\end{lemma}
\begin{lemma}\label{lemma4}
	Let $G(x)$ be an integrable function on $\mathbf{R}$ of bounded variation $||G||_V,\ X$ be a random variable and $a,b \not= 0$ are real numbers. Then \[\mathbf{E}\bigg[G\bigg(\frac{X+a}{b}\bigg)\bigg] \le ||G||_V\sup_u\Big|\mathbf{P}\big(X \le u\big) - \Phi(u)\Big| + ||G||_1|b|,\] where $||G||_1$ is the $L_1(\mathbf{R})$ norm of $G(x).$	
\end{lemma}

In the proof of Theorem \ref{th2}, we also need the following lemma of El Machkouri and Ouchti \cite{M07}.
\begin{lemma}\label{lemma5}
	Let $X$ and $Y$ be two random variables. Then, for $p\geq 1$,
	\begin{equation}
	D(X+Y)\leq 2D(X)+3 \Big\|\mathbf{E}\big[Y^{2p}|X\big]\Big\|_1^{1/{(2p+1)}}.
	\end{equation}
\end{lemma}

\subsection{Proof of Theorem \ref{th1}}\label{sec3}
By Lemma \ref{lemma1}, we only need to consider the case of $ \rho \in (0,1]. $
We follow the method of Grama and Haeusler \cite{GH}.
Let $T = 1 + \delta_n^2.$ We introduce a modification of the conditional variance  $\langle X \rangle_n$ as follows:
\begin{equation}
V_k = \langle X \rangle_k\mathbf{1}_{\{k< n\}} + T\mathbf{1}_{\{k=n\}}.
\end{equation}
It is easy to see that $V_0 = 0, V_n = T$, and that $(V_k,\mathcal{F}_k)_{k=0,...,n}$ is a predictable process. Set \[\gamma = \epsilon_n + \delta_n.\] Let $c_*$ be some positive and sufficient large constant. Define the following non-increasing discrete time predictable process
\begin{equation}
A_k = c_*^2\gamma^2 + T - V_k,\quad k = 1,...,n.
\end{equation}
Obviously, we have $A_0 = c_*^2\gamma^2 + T$ and $A_n = c_*^2\gamma^2.$ In addition, for $u,x \in \mathbf{R}$, and $y > 0$, denote
\begin{equation}
\Phi_u(x,y) = \Phi\bigg(\frac{u-x}{\sqrt{y}}\bigg).
\end{equation}

Let $\mathcal{N} = \mathcal{N}(0,1)$ be a standard normal random variable, which is independent of $X_n$. Using a smoothing procedure, by Lemma \ref{lemma3}, we deduce that
\begin{eqnarray}\label{great1}
\sup_u\Big|\mathbf{P}\big(X_n \le u\big) - \Phi(u)\Big| &\le& c_1\sup_{u}\Big|\mathbf{P}\big(X_n + c_*\gamma\mathcal{N} \le u\big) - \Phi(u)\Big| + c_2\gamma  \nonumber\\
&=&c_1\sup_u\Big|\mathbf{E}\big[\Phi_u\big(X_n,A_n\big)\big] - \Phi(u)\Big| + c_2\gamma \nonumber\\
&\le& c_1\sup_u\Big|\mathbf{E}\big[\Phi_u\big(X_n,A_n\big)\big] - \mathbf{E}\big[\Phi_u\big(X_0,A_0\big)\big]\Big| \nonumber\\
&& +c_1\sup_u\Big|\mathbf{E}\big[\Phi_u\big(X_0,A_0\big)\big] - \Phi(u)\Big| + c_2\gamma\nonumber\\
&=& c_1\sup_u\Big|\mathbf{E}\big[\Phi_u\big(X_n,A_n\big)\big] - \mathbf{E}\big[\Phi_u\big(X_0,A_0\big)\big]\Big| \nonumber\\
&& +c_1\sup_u\bigg|\Phi\bigg(\frac{u}{\sqrt{c_*^2\gamma^2 + T}}\bigg) - \Phi(u)\bigg| + c_2\gamma.
\end{eqnarray}
It is obvious that
\begin{equation}
\bigg|\Phi\bigg(\frac{u}{\sqrt{c_*^2\gamma^2 + T}}\bigg) - \Phi(u)\bigg| \le c_3\bigg|\frac{1}{\sqrt{c_*^2\gamma^2 + T}} - 1\bigg| \le c_4\gamma.
\end{equation}
Returning to (\ref{great1}), we get
\begin{equation}\label{great2}
\sup_u\Big|\mathbf{P}\big(X_n \le u\big) - \Phi(u)\Big| \le  c_1\sup_u\Big|\mathbf{E}\big[\Phi_u\big(X_n,A_n\big)\big] - \mathbf{E}\big[\Phi_u\big(X_0,A_0\big)\big]\Big| + c_5\gamma.
\end{equation}
By a simple telescoping, we know that
\begin{eqnarray}
\mathbf{E}\big[\Phi_u\big(X_n,A_n\big)\big] - \mathbf{E}\big[\Phi_u\big(X_0,A_0\big)\big] &= \mathbf{E}\Big[\sum_{k=1}^{n}\Big(\Phi_u\big(X_k,A_k\big) - \Phi_u\big(X_{k-1},A_{k-1}\big) \Big)\Big]. 		
\end{eqnarray}
Taking into account the fact that
\begin{eqnarray*}
\frac{\partial^2}{\partial x^2}\Phi_u(x,y) = 2\frac{\partial}{\partial y}\Phi_u(x,y),
\end{eqnarray*}
we get
\begin{equation}\label{great3}
\mathbf{E}\big[\Phi_u\big(X_n,A_n\big)\big] - \mathbf{E}\big[\Phi_u\big(X_0,A_0\big)\big] = J_1 + J_2 - J_3, 		
\end{equation}
where
\begin{eqnarray}
J_1 &=& \mathbf{E}\bigg[\sum_{k=1}^{n}\bigg(\Phi_u(X_k,A_k) - \Phi_u(X_{k-1},A_k) - \frac{\partial}{\partial x}\Phi_u(X_{k-1},A_k)\xi_k \nonumber\\&\;&  \quad\quad\quad\quad \quad\quad\quad\quad - \frac{1}{2}\frac{\partial^2}{\partial x^2}\Phi_u(X_{k-1},A_k)\xi_k^2 - \frac{1}{6}\frac{\partial^3}{\partial x^3}\Phi_u(X_{k-1},A_k)\xi_k^3\bigg)\bigg],\\
J_2 &=& \frac{1}{2}\mathbf{E}\bigg[\sum_{k=1}^{n}\frac{\partial^2}{\partial x^2}\Phi_u(X_{k-1},A_k)\Big(\bigtriangleup \langle X \rangle_k - \bigtriangleup V_k \Big)\bigg],\\
J_3 &=& \mathbf{E}\bigg[\sum_{k=1}^{n}\bigg(\Phi_u(X_{k-1},A_{k-1}) - \Phi_u(X_{k-1},A_k) - \frac{\partial}{\partial y}\Phi_u(X_{k-1},A_k)\bigtriangleup V_k \bigg)\bigg],	
\end{eqnarray}
where $ \bigtriangleup \langle X \rangle_k =\langle X \rangle_k - \langle X \rangle_{k-1}.  $

Now, we need to give some estimates of $J_1,J_2$ and $J_3$. To this end, we introduce some notations. Denote by $\vartheta_i$ some random variables  satisfying $0 \le \vartheta_i \le 1$, which may represent different values at different places. For the rest of the paper, $\varphi$ stands for the density function of the standard normal random variable.

\emph{Control of $ J_1 :$} For convenience's sake, let $T_{k-1} = \big(u - X_{k-1}\big)/\sqrt{A_{k}}, \  k = 1,2,...,n.$ It is easy to see that
\begin{eqnarray}
B_k &=:& \Phi_u(X_k,A_k) - \Phi_u(X_{k-1},A_k) - \frac{\partial}{\partial x}\Phi_u(X_{k-1},A_k)\xi_k\nonumber\\&\;& - \frac{1}{2}\frac{\partial^2}{\partial x^2}\Phi_u(X_{k-1},A_k)\xi_k^2 - \frac{1}{6}\frac{\partial^3}{\partial x^3}\Phi_u(X_{k-1},A_k)\xi_k^3\nonumber\\
&=& \Phi\bigg(T_{k-1} - \frac{\xi_k}{\sqrt{A_k}}\bigg) - \Phi(T_{k-1}) + \Phi'(T_{k-1})\frac{\xi_k}{\sqrt{A_k}}\nonumber\\&\;& - \frac{1}{2}\Phi''(T_{k-1})\bigg(\frac{\xi_k}{\sqrt{A_k}}\bigg)^2 + \frac{1}{6}\Phi'''(T_{k-1})\bigg(\frac{\xi_k}{\sqrt{A_k}}\bigg)^3.\nonumber
\end{eqnarray}
To estimate the right hand side of the last equality, we distinguish two cases.

\emph{Case 1: }$ |\xi_k/\sqrt{A_k}| \le 2 + |T_{k-1}|/2 .$ By a four-term Taylor expansion, it is obvious that if $ |\xi_k/\sqrt{A_k}| \le 1, $ then
\begin{eqnarray}
\bigg|B_k\bigg| &=& \bigg|\frac{1}{24}\Phi^{(4)}\bigg(T_{k-1} - \vartheta\frac{\xi_k}{\sqrt{A_k}}\bigg)\bigg|\frac{\xi_k}{\sqrt{A_k}}\bigg|^4\bigg|\nonumber\\
&\le&\bigg|\Phi^{(4)}\bigg(T_{k-1} - \vartheta\frac{\xi_k}{\sqrt{A_k}}\bigg)\bigg|\bigg|\frac{\xi_k}{\sqrt{A_k}}\bigg|^{3 + \rho}.\nonumber
\end{eqnarray}
If $ |\xi_k/\sqrt{A_k}|\; >\; 1, $ by a three-term Taylor expansion, then
\begin{eqnarray}
\bigg|B_k\bigg| &\le& \frac{1}{2}\bigg(\bigg|\Phi'''\bigg(T_{k-1} - \vartheta\frac{\xi_k}{\sqrt{A_k}}\bigg)\bigg| + \bigg|\Phi'''(T_{k-1})\bigg|\bigg)\bigg|\frac{\xi_k}{\sqrt{A_k}}\bigg|^3 \nonumber\\
&\le& \bigg|\Phi'''\bigg(T_{k-1} - \vartheta'\frac{\xi_k}{\sqrt{A_k}}\bigg)\bigg|\bigg|\frac{\xi_k}{\sqrt{A_k}}\bigg|^3 \nonumber\\&\le&\bigg|\Phi'''\bigg(T_{k-1} - \vartheta'\frac{\xi_k}{\sqrt{A_k}}\bigg)\bigg|\bigg|\frac{\xi_k}{\sqrt{A_k}}\bigg|^{3+\rho}, \nonumber		
\end{eqnarray}
where
\begin{displaymath}
 \vartheta' =  \left\{ \begin{array}{ll}
\vartheta,  & \textrm{if $ \big|\Phi'''\big(T_{k-1} - \vartheta\frac{\xi_k}{\sqrt{A_k}}\big)\big| \ge |\Phi'''(T_{k-1})|,$}\\
  0, & \textrm{if $\big|\Phi'''\big(T_{k-1} - \vartheta\frac{\xi_k}{\sqrt{A_k}}\big)\big|   <  |\Phi'''(T_{k-1})|$}.
\end{array} \right.
\end{displaymath}
Using the inequality $ \max\{|\Phi'''(t)|,|\Phi''''(t)|\} \le \varphi(t)(2 + t^4), $ we find that
\begin{eqnarray}\label{B1}
\Big|B_k\mathbf{1}_{\{|\xi_k/\sqrt{A_k}| \le 2 + |T_{k-1}|/2\}}\Big| &\le& \varphi\bigg(T_{k-1} - \vartheta_1\frac{\xi_k}{\sqrt{A_k}}\bigg)\Bigg(2 + \bigg(T_{k-1} - \vartheta_1\frac{\xi_k}{\sqrt{A_k}}\bigg)^4\Bigg)\Big|\frac{\xi_k}{\sqrt{A_k}}\Big|^{3+\rho}\nonumber \\&\le& g_1(T_{k-1})\Big|\frac{\xi_k}{\sqrt{A_k}}\Big|^{3+\rho},
\end{eqnarray}
where\[ g_1(z) = \sup_{|t-z|\le 2 + |z|/2}\varphi(t)(2+t^4). \]

\emph{Case 2: }$ |\xi_k/\sqrt{A_k}| > 2 + |T_{k-1}|/2 .$ It is obvious that, for $ |\bigtriangleup x| > 1 + |x|/2, $
\begin{eqnarray*}
&&\Big|\Phi(x - \bigtriangleup x) - \Phi(x) + \Phi'(x)\bigtriangleup x - \frac{1}{2}\Phi''(x)(\bigtriangleup x)^2 + \frac{1}{6}\Phi'''(x)(\bigtriangleup x)^3\Big| \\&&\le \bigg(\bigg|\frac{\Phi(x - \bigtriangleup x) - \Phi(x)}{|\bigtriangleup x|^{3}}\bigg| + |\Phi'(x)| + |\Phi''(x)| + |\Phi'''(x)|\bigg)|\bigtriangleup x|^{3}\\
&&\le \bigg(8\bigg|\frac{\Phi(x - \bigtriangleup x) - \Phi(x)}{(2 + |x|)^{3}}\bigg|+ |\Phi'(x)| + |\Phi''(x)| + |\Phi'''(x)|\bigg)|\bigtriangleup x|^{3}\\
&&\le \bigg(\frac{\widetilde{c}}{(2 + |x|)^{3}} + |\Phi'(x)| + |\Phi''(x)| + |\Phi'''(x)|\bigg)|\bigtriangleup x|^{3}\\
&&\le\frac{\hat{c}}{(2 + |x|)^{3}}|\bigtriangleup x|^{3}\\
&&\le \frac{\hat{c}}{(2 + |x|)^{3}}|\bigtriangleup x|^{3 + \rho}.
\end{eqnarray*}
Hence, we have
\begin{eqnarray}\label{B2}
\Big|B_k\mathbf{1}_{\{|\xi_k/\sqrt{A_k}| > 2 + |T_{k-1}|/2\}}\Big|  \le g_2(T_{k-1})\Big|\frac{\xi_k}{\sqrt{A_k}}\Big|^{3+\rho},
\end{eqnarray}
where
\[ g_2(z) = \frac{\hat{c}}{(2 + |z|)^{3}}. \]

Denote $${G}(z) = g_1(z) + g_2(z).$$
Combining (\ref{B1}) and (\ref{B2}) together, we get
\begin{equation}\label{J11}
|B_k| \le G(T_{k-1})\bigg|\frac{\xi_k}{\sqrt{A_k}}\bigg|^{3 + \rho}.
\end{equation}
Therefore,
\begin{equation}\label{J113}
\bigg|J_1\bigg| = \bigg|\mathbf{E}\bigg[\sum_{k=1}^{n}B_k\bigg]\bigg| \le \mathbf{E}\bigg[\sum_{k=1}^{n}G(T_{k-1})\bigg|\frac{\xi_k}{\sqrt{A_k}}\bigg|^{3 + \rho}\bigg].
\end{equation}
Next, we consider conditional expectation of $|\xi_k|^{3 + \rho}.$ By condition (\ref{condi3}), we get
\begin{equation}
\mathbf{E}\big[|\xi_k|^{3+\rho}\big|\mathcal{F}_{k-1}\big] \le \epsilon_n^{1+\rho}\bigtriangleup\langle X \rangle_k,
\end{equation}
where $\bigtriangleup\langle X \rangle_k = \langle X \rangle_k - \langle X \rangle_{k-1}.$ And we know that
\begin{equation}\label{J12}
\bigtriangleup\langle X \rangle_k = \bigtriangleup V_k = V_k - V_{k-1}, \ 1\le k < n,\ \bigtriangleup\langle X \rangle_n \le \bigtriangleup V_n,
\end{equation}
then
\begin{equation}\label{J13}
\mathbf{E}\big[|\xi_k|^{3+\rho}\big|\mathcal{F}_{k-1}\big] \le \epsilon_n^{1+\rho}\bigtriangleup V_k.
\end{equation}
By (\ref{J113}) and (\ref{J13}), we obtain
\begin{equation}\label{R1}
|J_1| \le R_1 := \epsilon_n^{1+\rho}\bigg[\sum_{k=1}^{n}\frac{G(T_{k-1})}{A_k^{(3+\rho)/2}}\bigtriangleup V_k\bigg].
\end{equation}

To estimate $ R_1 $, we introduce the time change $\tau_t$ as follow: for any real $t \in [0,T],$
\begin{equation}
\tau_t = \min\{k\le n: V_k \ge t\}, \ \  \textrm{where} \ \min\emptyset = n.
\end{equation}
Obviously, for any $t \in [0,T],$ the stopping time $\tau_t$ is predictable. In addition, $(\sigma_{k})_{k=1,...,n+1}(\textrm{with}\ \sigma_{1} = 0)$ stands for the increasing sequence of moments when the increasing and stepwise function $\tau_t,t \in [0,T]$, has jumps. It is easy to see that $\bigtriangleup V_k = \int_{[\sigma_{k},\sigma_{k+1})}dt$, and that $k = \tau_t$ for $t \in [\sigma_{k},\sigma_{k+1})$. Since $\tau_T = n$, we have
\begin{equation}
\sum_{k=1}^{n}\frac{G(T_{k-1})}{A_k^{(3 + \rho)/2}}\bigtriangleup V_k = \sum_{k=1}^{n}\int_{[\sigma_{k},\sigma_{k+1})}\frac{G(T_{\tau_t-1})}{A_{\tau_t}^{(3 + \rho)/2}}dt = \int_{0}^{T}\frac{G(T_{\tau_t-1})}{A_{\tau_t}^{(3 + \rho)/2}}dt.
\end{equation}
Let $a_t = c_*^2\gamma^2 + T - t$. Because of $\bigtriangleup V_{\tau_t} \le 2\epsilon_n^2 + 2\delta_n^2$ (cf. Lemma \ref{lemma2}), we know that
\begin{equation}
t \le  V_{\tau_t} = V_{\tau_t-1} + \bigtriangleup V_{\tau_t} \le t + 2\epsilon_n^2 + 2\delta_n^2,\quad  t \in [0,T].
\end{equation}
Assume   $c_* \ge 2$, then we have
\begin{equation}
\frac{1}{2}a_t \le A_{\tau_t} = c_*^2\gamma^2 + T - V_{\tau_t} \le a_t,\quad t \in [0,T].
\end{equation}
Note that ${G}(z)$ is symmetric and is non-increasing in $z \ge 0$. The last bound implies that
\begin{equation}\label{R2}
R_{1} \le 2^{(3 + \rho)/2}\epsilon_n^{1+\rho}\int_{0}^{T}\frac{1}{a_{t}^{(3 + \rho)/2}}\mathbf{E}\Big[G\Big(\frac{u - X_{\tau_t - 1}}{a_t^{1/2}}\Big)\Big]dt.
\end{equation}
Note also that $G(z)$ is a symmetric integrable function of bounded variation. By Lemma \ref{lemma4}, it is obvious that
\begin{equation}\label{R3}
\mathbf{E}\Big[G\Big(\frac{u - X_{\tau_t - 1}}{a_t^{1/2}}\Big)\Big] \le c_6\sup_z\Big|\mathbf{P}\big(X_{\tau_t-1} \le z\big) - \Phi(z)\Big| + c_7\sqrt{a_t}.
\end{equation}
Because of $c_* \ge 2, V_{\tau_t - 1} = V_{\tau_t} - \bigtriangleup V_{\tau_t},\ V_{\tau_t} \ge t$ and $\bigtriangleup V_{\tau_t} \le 2\epsilon_n^2 + 2\delta_n^2$, we obtain
\begin{equation}
V_n - V_{\tau_t-1} = V_n - V_{\tau_t} + \bigtriangleup V_{\tau_t} \le 2\epsilon_n^2 + 2\delta_n^2 + T - t \le a_t.
\end{equation}
Therefore
\begin{eqnarray*}
\mathbf{E}\Big[\big(X_n - X_{\tau_t-1}\big)^2\Big|\mathcal{F}_{\tau_t-1}\Big] &=& \mathbf{E}\bigg[\sum_{k=\tau_t}^{n}\mathbf{E}\big[\xi_k^2\big|\mathcal{F}_{k-1}\big]\bigg|\mathcal{F}_{\tau_t-1}\bigg]\\
&=&  \mathbf{E}\big[\langle X \rangle_n - \langle X \rangle_{\tau_t-1}\big|\mathcal{F}_{\tau_t-1}\big]\\
&\le& \mathbf{E}[ V_n - V_{\tau_t-1}|\mathcal{F}_{\tau_t-1}]\\
&\le& a_t.
\end{eqnarray*}
Then, by Lemma \ref{lemma3}, we deduce that for any $t \in [0,T]$,
\begin{equation}\label{R4}
\sup_z\big|\mathbf{P}\big(X_{\tau_t-1} \le z\big) - \Phi(z)\big| \le c_8\sup_z\big|\mathbf{P}\big(X_n \le z\big) - \Phi(z)\big| + c_9\sqrt{a_t}.
\end{equation}
Combining (\ref{R1}), (\ref{R2}), (\ref{R3})\ and\ (\ref{R4}) together, we get
\begin{equation}
|J_1| \le c_{10}\epsilon_n^{1 + \rho}\int_{0}^{T}\frac{1}{a_{t}^{(3 + \rho)/2}}dt\sup_z\big|\mathbf{P}\big(X_n \le z\big) - \Phi(z)\big| + c_{11}\epsilon_n^{1 + \rho}\int_{0}^{T}\frac{1}{a_t^{1 + \rho/2}}dt.
\end{equation}
Taking some elementary computations, it follows that
\begin{eqnarray}
\int_{0}^{T}\frac{1}{a_{t}^{(3 + \rho)/2}} dt = \int_{0}^{T}\frac{1}{(c_*^2\gamma^2 + T - t)^{(3 + \rho)/2}}dt \le \frac{2}{c_*^{1 + \rho}(1+\rho)\gamma^{1 + \rho}}
\end{eqnarray}
and
\begin{eqnarray}
  \int_{0}^{T}\frac{1}{a_{t}^{1 + \rho/2}}dt  = \int_{0}^{T}\frac{1}{(c_*^2\gamma^2 + T - t)^{1 + \rho/2}} dt\le \frac{2}{c_*^{\rho}\rho\gamma^{\rho}}.
\end{eqnarray}
This yields
\begin{equation}\label{J21}
\big|J_1\big| \le \frac{c_{12}}{c_*^{1+\rho}}\sup_z\big|\mathbf{P}\big(X_n \le z\big) - \Phi(z)\big| + \frac{c_{\rho,1}\epsilon_n}{\rho}.
\end{equation}

\emph{Control of $ J_2 :$} Since $ 0 \le \bigtriangleup V_k - \bigtriangleup \langle X \rangle_k \le 2\delta^2\mathbf{1}_{\{k=n\}}, $ we have\[ |J_2| \le \mathbf{E}\Big[\frac{1}{2A_n}\big|\varphi'(T_{n-1})(\bigtriangleup V_n - \bigtriangleup \langle X \rangle_n )\big|\Big]. \]
Denote $ \widetilde{G}(z) = \sup_{|z-t|\le 1}|\varphi'(t)|, $ and then $ |\varphi'(z)| \le \widetilde{G}(z) $ for any real $ z. $ Since $ A_n = c_*^2\gamma^2, $ then we get the following estimation:\[ |J_2 | \le \frac{1}{c_*^2}\mathbf{E}\big[\widetilde{G}(T_{n-1})\big]. \]
Note that $ \widetilde{G} $ is non-increasing in $ z \ge 0, $ and thus it has bounded variation on $\mathbf{R}.$ By Lemma \ref{lemma4}, we get
\begin{equation}
|J_2| \le \frac{c_{13}}{c_*^2}\sup_z\big|\mathbf{P}\big(X_{n-1} \le z\big) - \Phi(z)\big| + c_{*,2}(\epsilon_n + \delta_n).
\end{equation}
Then, by Lemma \ref{lemma3}, we deduce that
\begin{equation}
\sup_z\big|\mathbf{P}\big(X_{n-1} \le z\big) - \Phi(z)\big| \le c_{14}\sup_z\big|\mathbf{P}\big(X_{n} \le z\big) - \Phi(z)\big| + c_{15}\epsilon_n.
\end{equation}
This yields
\begin{equation}\label{J22}
|J_2| \le \frac{c_{16}}{c_*^2}\sup_z\big|\mathbf{P}\big(X_{n} \le z\big) - \Phi(z)\big| + c_{\rho,2}(\epsilon_n + \delta_n).
\end{equation}

\emph{Control of $ J_3 $.} By a two-term Taylor expansion, it follows that\[ |J_3| = \frac{1}{8}\mathbf{E}\bigg[\sum_{k=1}^{n}\frac{1}{(A_k - \vartheta_k \bigtriangleup A_k)^2}\varphi'''\bigg(\frac{u-X_{k-1}}{\sqrt{A_k - \vartheta_k \bigtriangleup A_k}}\bigg)(\bigtriangleup A_k)^2\bigg]. \]
Note that $ c_* \ge 2, \bigtriangleup A_k \le 0 $ and, by Lemma (\ref{lemma2}), $ |\bigtriangleup A_k| = \bigtriangleup V_k \le 2\epsilon_n^2 + 2\delta_n^2. $ We obtain
\begin{equation}\label{J}
A_k \le A_k - \vartheta_k \bigtriangleup A_k \le c_*^2\gamma^2 + T - V_k + 2\epsilon_n^2 + 2\delta_n^2 \le 2A_k.
\end{equation}
Denote $\widehat{G}(z) = \sup_{|t-z| \le 2}|\varphi'''(t)|$. Then $\widehat{G}(z)$ is symmetric, and is non-increasing in $z \ge 0$. Using (\ref{J}), we get
\begin{equation}
|J_3| \le (2\epsilon_n^2 + 2\delta_n^2) \mathbf{E}\Bigg[\sum_{k=1}^n\frac{1}{A_{k}^2}\widehat{G}\bigg(\frac{T_{k-1}}{\sqrt{2}}\bigg)\bigtriangleup V_k\Bigg].
\end{equation}
By  an argument  similar to that of (\ref{J21}), we get
\begin{eqnarray}\label{J23}
|J_3| &\le& \frac{c_{17}(2\epsilon_n^2 + 2\delta_n^2)}{c_*^2\gamma^2}\sup_z\big|\mathbf{P}\big(X_n \le z\big) - \Phi(z)\big| + \frac{2c_{18}(2\epsilon_n^2 + 2\delta_n^2)}{c_*\gamma} \nonumber\\ &\le& \frac{c_{19}}{c_*^2}\sup_z\big|\mathbf{P}\big(X_n \le z\big) - \Phi(z)\big| + \frac{4c_{18}(\epsilon_n + \delta_n)^2}{c_*\gamma} \nonumber \\&\le& \frac{c_{19}}{c_*^2}\sup_z\big|\mathbf{P}\big(X_n \le z\big) - \Phi(z)\big| + c_{\rho,3}(\epsilon_n + \delta_n).
\end{eqnarray}

Combining (\ref{great3}), (\ref{J21}), (\ref{J22}) and (\ref{J23}) together, we get
\[
\Big|\mathbf{E}\big[\Phi_u\big(X_n,A_n\big)\big] - \mathbf{E}\big[\Phi_u\big(X_0,A_0\big)\big]\Big| \le \frac{c_{20}}{c_*^{1+\rho}}\sup_z\big|\mathbf{P}\big(X_n \le z\big) - \Phi(z)\big| + \frac{\hat{c}_{\rho}}{\rho}(\epsilon_n + \delta_n),
\]
By (\ref{great2}), we know that
\[
\sup_z\big|\mathbf{P}\big(X_n \le z\big) - \Phi(z)\big| \le \frac{c_{21}}{c_*^{1+\rho}}\sup_z\big|\mathbf{P}\big(X_n \le z\big) - \Phi(z)\big| + \frac{\tilde{c}_{\rho}}{\rho}(\epsilon_n + \delta_n),
\]
from which, choosing $c_*^{1+\rho} = \max{\{2c_{21},2^{1+\rho}}\},$ we get
\begin{equation}
\sup_z\big|\mathbf{P}\big(X_n \le z\big) - \Phi(z)\big| \le \frac{2\tilde{c}_{\rho}(\epsilon_n + \delta_n)}{\rho}.
\end{equation}

\subsection{Proof of Theorem \ref{th2}}
Following the method of Bolthausen \cite{B82}, we enlarge the sequence $(\xi_i, \mathcal{F}_i)_{1\leq i\leq n}$ to $\big(\hat{\xi}_i, \hat{\mathcal{F}}_{i}\big)_{1\leq i\leq N}$ such that $\big\langle\hat{X}\big\rangle_N:= \sum_{i=1}^{N}\mathbf{E}\big[\hat{\xi}_i^2|\hat{\mathcal{F}}_{i-1}\big] = 1$ a.s., and then apply Theorem \ref{th1} to the enlarged sequence.
Consider the stopping time
\begin{equation}
\tau=\sup\{k\leq n: \langle X\rangle_k\leq 1\}.
\end{equation}
Assume that $0\leq\varepsilon\leq\epsilon_n.$ Let $r=\Big\lfloor\frac{1-\langle X\rangle_\tau}{\varepsilon^2}\Big\rfloor$, where $\lfloor x\rfloor$ denotes the ``integer part" of $x$. It is easy to see that $r\leq\Big\lfloor \frac{1}{\varepsilon^2}\Big\rfloor.$
Set $N=n+r+1.$ Let $ (\zeta_i)_{i \ge 1} $ be a sequence of independent Rademacher random variables, which is   independent of the martingale differences $ (\xi_i)_{1 \le i \le n}. $
Consider the random variables  $\big(\hat{\xi}_i, \hat{\mathcal{F}}_{i}\big)_{1\leq i\leq N} $ defined as follows:
\begin{displaymath}
\hat{\xi_i} =  \left\{ \begin{array}{ll}
\xi_i \ \textrm{\ \ a.s.},& \textrm{if}\ i \le \tau, \\
\varepsilon\zeta_i \ \textrm{\ \ a.s.},& \textrm{if}\ \tau+1\leq i\leq \tau +r,   \\
\big(1- \langle X \rangle_\tau -r\varepsilon^2\big)^{1/2}\zeta_i \ \textrm{\ \ a.s.},& \textrm{if}\ i = \tau + r +1, \\
0 \ \textrm{\ \ a.s.},& \textrm{if}\ \tau + r + 1 \leq i \leq N,
\end{array} \right.
\end{displaymath}
and $ \hat{\mathcal{F}}_i = \sigma\big(\hat{\xi}_1, \hat{\xi}_2,...,\hat{\xi}_i\big). $

Clearly, $\big(\hat{\xi}_i, \hat{\mathcal{F}}_i\big)_{1\leq i\leq N}$ still forms a martingale difference sequence with respect to the enlarged filtration. Then $ \hat{X}_k = \sum_{i=1}^{k}\hat{\xi}_i,\ k = 0,...,N $, with $ \hat{X}_0 = 0 $, is also a martingale. Moreover, it holds that $\big\langle\hat{X}\big\rangle_N=1$, $\mathbf{E}\big[\hat{\xi}_i^3\big|\hat{\mathcal{F}}_{i-1}\big]=0$ and
$$\mathbf{E}\big[\big|\hat{\xi}_i\big|^{3+\rho}\big|\hat{\mathcal{F}}_{i-1}\big]\leq \epsilon_n^{1+\rho}\mathbf{E}\big[\hat{\xi}_i^2\big|\hat{\mathcal{F}}_{i-1}\big], \quad \ \ a.s.$$
By Theorem \ref{th1}, we have
\begin{equation}
D\big(\hat{X}_N\big)\leq \frac{c_{\rho}\epsilon_n}{\rho}.
\end{equation}
Using Lemma \ref{lemma5}, we obtain that
\begin{eqnarray}\label{2}
D(X_n)
& \leq& 2D\big(\hat{X}_N\big)+3\Big\|\mathbf{E}\Big[\Big|X_n-\hat{X}_N\Big|^{2p}\Big|\hat{X}_N\Big]\Big\|_{1}^{1/{(2p+1)}}\nonumber\\
& \leq&  \frac{2c_{\rho}\epsilon_n}{\rho} + 3\Big(\mathbf{E}\Big[\Big|\hat{X}_N-X_n\Big|^{2p}\Big]\Big)^{1/{(2p+1)}}.
\end{eqnarray}
Since $\tau$ is a stopping time and
\begin{equation}
\hat{X}_N-X_n=\sum_{i=\tau+1}^N\Big(\hat{\xi}_i-\xi_i\Big), \qquad \mbox{ where put } \xi_i=0 \mbox{ for } i>n,
\end{equation}
$(\hat{\xi}_i-\xi_i, \hat{\mathcal{F}}_i)_{i\geq \tau+1}$ still forms a martingale difference sequence. Applying Burkhold's inequality (cf. Theorem 2.11 of Hall and Heyde \cite{H80}), we get
\begin{eqnarray}\label{E1}
\mathbf{E}\Big[\Big|\hat{X}_N-X_n\Big|^{2p}\Big]&\le& \mathbf{E}\Big[\max_{\tau+1\leq i\leq N}\Big|\hat{X}_i-X_i\Big|^{2p}\Big]\nonumber \\&\leq& c_p\Bigg(\mathbf{E}\Big[\Big|\sum_{i=\tau+1}^N\mathbf{E}\Big[\Big(\hat{\xi}_i-\xi_i\Big)^2\Big|\hat{\mathcal{F}}_{i-1}\Big]\Big|^p\ \Big]+\mathbf{E}\Big[\max_{\tau+1\leq i\leq N}\Big|\hat{\xi}_i-\xi_i\Big|^{2p}\Big]\Bigg).
\end{eqnarray}
As $\xi_i$ and $\hat{\xi}_i$ be orthogonal random variables, we have
\begin{eqnarray*}
\sum_{i=\tau+1}^N\mathbf{E}\bigg[\Big(\hat{\xi}_i-\xi_i\Big)^2\bigg|\hat{\mathcal{F}}_{i-1}\bigg]=\sum_{i=\tau+1}^N\mathbf{E}\Big[\hat{\xi}_i^2\Big|\hat{\mathcal{F}}_{i-1}\Big]+\sum_{i=\tau+1}^n\mathbf{E}\Big[\xi_i^2\Big|\hat{\mathcal{F}}_{i-1}\Big]=1-2\langle X\rangle_\tau+\langle X\rangle_n.
\end{eqnarray*}
Noting that $1-\mathbf{E}[\xi_{\tau+1}^2|\mathcal{F}_{\tau}] \le \langle X\rangle_\tau.$ Consequently, using the inequality $|a+b|^{p} \leq 2^{p-1}\left(|a|^{p}+|b|^{p}\right), p\geq1$, and  Jensen's inequality, we derive that
\begin{eqnarray}
\Bigg|\sum_{i=\tau+1}^N\mathbf{E}\Big[\Big(\hat{\xi}_i-\xi_i\Big)^2\Big|\hat{\mathcal{F}}_{i-1}\Big]\Bigg|^p &\leq& \Big|\langle X\rangle_n-1+2\mathbf{E}\big[\xi_{\tau+1}^2\big|\mathcal{F}_{\tau}\big]\Big|^p\nonumber\\
&\leq& 2^{2p-1}\bigg(\big|\langle X\rangle_n-1\big|^p+\Big|\mathbf{E}\big[\xi_{\tau+1}^2\big|\mathcal{F}_{\tau}\big]\Big|^p\bigg)\nonumber\\
&\leq& 2^{2p-1}\Big(\big|\langle X\rangle_n-1\big|^p+\mathbf{E}\big[|\xi_{\tau+1}|^{2p} \big|\mathcal{F}_{\tau}\big]\Big).
\end{eqnarray}
Taking expectations on both sides of the last inequality, we deduce that
\begin{eqnarray}\label{E11}
\mathbf{E}\Bigg[\Bigg|\sum_{i=\tau+1}^N\mathbf{E}\bigg[\Big(\hat{\xi}_i-\xi_i\Big)^2\bigg|\hat{\mathcal
{F}}_{i-1}\bigg]\Bigg|^p\ \Bigg]
&\leq& 2^{2p-1} \bigg(\mathbf{E}\Big[\big|\langle X\rangle_n-1\big|^p\Big]+\mathbf{E}\big[|\xi_{\tau+1}|^{2p}\big]\bigg)\nonumber\\
&\leq& 2^{2p-1} \bigg(\mathbf{E}\Big[\big|\langle X\rangle_n-1\big|^p\Big]+\mathbf{E}\Big[\max_{1\leq i\leq n}|\xi_i|^{2p}\Big]\bigg).
\end{eqnarray}
Similarly, using the inequality $|a+b|^{p} \leq 2^{p-1}\left(|a|^{p}+|b|^{p}\right), p\geq1$,
\begin{eqnarray}\label{E12}
\mathbf{E}\bigg[\max_{\tau+1\leq i\leq N}\Big|\hat{\xi}_i-\xi_i\Big|^{2p}\bigg]&\leq& 2^{2p-1} \mathbf{E}\bigg[\max_{\tau+1\leq i\leq N}\Big(|\xi_i|^{2p}+\big|\hat{\xi_i}\big|^{2p}\Big)\bigg]\nonumber\\
&\leq& 2^{2p-1} \bigg(\mathbf{E}\Big[\max_{1\leq i\leq n}|\xi_i|^{2p}\Big]+\varepsilon^{2p}\bigg).
\end{eqnarray}
Combining (\ref{E1}), (\ref{E11}) and (\ref{E12}) together, we obtain
\begin{equation}
\mathbf{E}\bigg[\Big|\hat{X}_N-X_n\Big|^{2p}\bigg] \leq \hat{c}_p\bigg(\mathbf{E}\Big[\big|\langle X\rangle_n-1\big|^{p}\Big] + \mathbf{E}\Big[\max_{1\leq i\leq n}|\xi_i|^{2p}\Big]+\varepsilon^{2p}\bigg).
\end{equation}
Finally, applying the last inequality to  (\ref{2})  and let  $\varepsilon\to 0$, then we have
\[ D(X_n)\leq \tilde{c}_{\rho}\frac{\epsilon_n}{\rho}+\tilde{c}_p\bigg(\mathbf{E}\Big[\big|\langle X\rangle_n - 1\big|^p\Big]+\mathbf{E}\Big[\max_{1\leq i\leq n}|\xi_i|^{2p}\Big]\bigg)^{1/{(2p+1)}}.
 \]
This completes the proof of  Theorem \ref{th2}.
\section*{Acknowledgements}%
The work has been partially supported by the National Natural Science Foundation
of China (Grant no.\,11601375). The research of Hailin Sang is partial supported by
the Simons Foundation (Grant no.\,586789).

\selectlanguage{english}

\end{document}